\documentclass{article}[12pt]
\usepackage{amssymb}
\usepackage{eurosym}
\usepackage{epsfig}
\usepackage{amsmath}
\usepackage{amsfonts}
\usepackage{pgf,tikz}
\usepackage{enumerate}
\usepackage{cite}
\usepackage[colorlinks=true, linkcolor=black, citecolor=black, urlcolor=black]{hyperref}

\newcommand{\R}{{{\Bbb R}}}

\newcommand{\N}{{{\Bbb N}}}

\newtheorem{theorem}{\sc Theorem}[section]
\newtheorem{proposition}{\sc Proposition}[section]
\newtheorem{lemma}{\sc Lemma}[section]
\newtheorem{definition}{\sc Definition}[section]
\newtheorem{remark}{\sc Remark}[section]
\newtheorem{corollary}{\sc Corollary}[section]
\newtheorem{example}{\sc Example}[section]

\def\qed{\hbox to 0pt{}\hfill$\rlap{$\sqcap$}\sqcup$\medbreak}

\title{On discontinuous differential equations and the method of solution-regions}

\author{Jorge Rodr\'iguez--L\'opez} 

\date{}
\begin{document}
 \maketitle

\begin{center}  {\small Departamento de Estat\'istica, An\'alise Matem\'atica e Optimizaci\'on, \\ Instituto de Matem\'aticas, Universidade de Santiago de Compostela, \\ 15782, Facultade de Matem\'aticas, Campus Vida, Santiago, Spain.\\  Email: jorgerodriguez.lopez@usc.es}
\end{center}

\medbreak

\noindent {\it Abstract.} We adapt the method of solution regions to prove new existence and localization results for systems of discontinuous differential equations. Some assumptions concerning the definition of a solution region are relaxed and thus our results enlarge the applicability of this method even in the case of Carath\'eodory nonlinearities. Several examples are provided to illustrate the theory.

\medbreak


\medbreak


\section{Introduction}

This paper deals with the existence of absolutely continuous solutions to the initial value problem
\begin{equation}\label{ivp}
x'(t)=f(t,x(t)) \quad \text{for a.a. } t\in I=[0,T], \quad x(0)=x_0\in\mathbb{R}^n,
\end{equation}
where $T>0$ and $n\in\N$ are fixed, and the function $f:I\times\R^n\rightarrow\R^n$ need not be continuous. 

It is well-known that if $f$ is not continuous with respect to the second variable, problem \eqref{ivp} may not have absolutely continuous solutions. Hence, we need to impose some assumption on $f$ at its discontinuity points. This assumption was inspired by a global transversality condition due to Bressan and Shen \cite{BreSh}, recently improved in \cite{pr}. 

Here, this technique is combined with the method of solution-regions due to Frigon \cite{f} (see also \cite{ftt,t}), which allow us to localize the solutions in a given compact set $R\subset I\times\R^n$. The set $R$ will be called a solution-region of problem \eqref{ivp}. As particular cases, it could be the set between ordered lower and upper solutions or the set given by a solution-tube, see \cite{f}. In this way, no monotonicity or growth conditions are imposed on $f$ in order to obtain the existence results. 

Our results enlarge the applicability of the method of solution-regions even in the case of initial value problems with continuous nonlinearities. We notice that one basic assumption in the notions of admissible and weak admissible pairs necessary to define the solution-region can be completely removed (compare Definitions \ref{def_adm} and \ref{def_Wadm} with Definition \ref{def_viable} below).   

The paper is organized as follows. In Section \ref{sec_pre} we present the notion of viable region and we recall some results about differential inclusions, which we employ in order to establish our existence result for \eqref{ivp}. In Section \ref{sec_exis} we state and prove the main result in the paper: an existence result for the initial value problem \eqref{ivp} which, in addition, provides localization in a solution-region. Finally, in Section \ref{sec_examples}, the applicability of the theory is illustrated by several meaningful examples. We focus on the localization of the solutions in a ball and the method of lower and upper solutions, even allowing them to be piecewise continuous functions.

\section{Preliminaries}\label{sec_pre}

\subsection{Admissible, weak admissible and viable regions}

First, let us introduce the notion of a Carath\'eodory function.

\begin{definition}
	A map $f:D\subset I\times\R^n\rightarrow\R^m$ is a Carath\'eodory function if 
	\begin{enumerate}
		\item[$(i)$] $f(t,\cdot)$ is continuous on $D_t=\{x:(t,x)\in D \}$ for almost every $t\in I$;
		\item[$(ii)$] $f(\cdot,x)$ is measurable for all $x\in \bigcup_{t\in I}D_t$;
		\item[$(iii)$] for all $k>0$, there exists $\psi_k\in L^1(I,\R)$ such that $\left\|f(t,x)\right\|\leq \psi_k(t)$ for a.e. $t$ and every $x$ such that $\left\|x\right\|\leq k$ and $(t,x)\in D$.
	\end{enumerate}
	A map $f:D\rightarrow\R^m$ is a locally Carath\'eodory function if $f:A\rightarrow\R^m$ is a Carath\'eodory function for every compact set $A\subset D$.
\end{definition}

Now, we recall the notions of admissible and weak admissible regions introduced in \cite{f} and \cite{t}, respectively.

\begin{definition}\label{def_adm}
	We say that a set $R\subset I\times\R^n$ is an \textit{admissible region} if there exist two continuous maps $h:I\times\R^n\rightarrow\R$ and $p=(p_1,p_2):I\times\R^n\rightarrow I\times\R^n$ satisfying the following conditions:
	\begin{enumerate}
		\item[$(i)$] $R=\{(t,x):h(t,x)\leq 0 \}$ is bounded and, for every $t\in I$,
		\[R_t=\{x\in\R^n:(t,x)\in R \}\neq \emptyset; \] 
		\item[$(ii)$]  the map $h$ has partial derivatives at $(t,x)$ for almost every $t$ and every $x$ with $(t,x)\in R^{c}:=(I\times \R^n)\setminus R$, and $\frac{\partial h}{\partial t}$ and $\nabla_x h$ are locally Carath\'eodory maps on $R^c$;
		\item[$(iii)$] $p$ is bounded and such that $p(t,x)=(t,x)$ for every $(t,x)\in R$ and 
		\[\langle \nabla_x h(t,x),p_2(t,x)-x \rangle<0 \quad \text{ for a.e. } t \text{ and every } x \text{ with } (t,x)\in R^c. \]
	\end{enumerate}
\end{definition}

\begin{definition}\label{def_Wadm}
	We say that a set $R\subset I\times\R^n$ is a \textit{weak admissible region} if there exist two continuous maps $h:I\times\R^n\rightarrow\R$ and $p=(p_1,p_2):I\times\R^n\rightarrow I\times\R^n$ satisfying conditions $(i)$ and $(ii)$ of Definition \ref{def_adm} and 
	\begin{enumerate}
		\item[$(iii^*)$] $p$ is bounded and such that $p(t,x)=(t,x)$ for every $(t,x)\in R$ and 
		\[\langle \nabla_x h(t,x),p_2(t,x)-x \rangle\leq 0 \quad \text{ for a.e. } t \text{ and every } x \text{ with } (t,x)\in R^c. \]
	\end{enumerate}
\end{definition}

Basically the difference between both definitions is the inequality in conditions $(iii)$ and $(iii^*)$, which can be non-strict in the weak notion. 

We will show here that these inequalities can be completely removed. In order to do that, we introduce the following concept.

\begin{definition}\label{def_viable}
	We say that  a set $R\subset I\times\R^n$ is a \textit{viable region} if there exist two continuous maps $h:I\times\R^n\rightarrow\R$ and $p:I\times\R^n\rightarrow I\times\R^n$, $p(t,x)=(p_1(t),p_2(t,x))$, satisfying conditions $(i)$ and $(ii)$ of Definition \ref{def_adm} and 	
	\begin{enumerate}
		\item[$(iii^\prime)$] $p$ is bounded and such that $p(t,x)=(t,x)$ for every $(t,x)\in R$.
	\end{enumerate}
	We call $(h,p)$ a \textit{viable pair} associated to $R$.
\end{definition}

It was shown in \cite{t} that any compact set $R\subset I\times\R^n$ such that $R_t\neq\emptyset$ for every $t\in I$ is a weak admissible region, so the definition of a viable region is very general. Notice that for some meaningful particular cases, the map $p$ can be seen as a retraction or as a projection of $I\times\R^n$ onto $R$. On the other hand, if $h$ is taken as a nonnegative function, then it could measure some distance between the point $(t,x)$ and the set $R$.

\subsection{Krasovskij solutions}

Associated to the initial value problem \eqref{ivp}, we consider the differential inclusion
\begin{equation}
\label{ivpI}
x'(t)\in\mathcal{K}f(t,x(t)) \quad \text{for a.a. } t \in I, \ \ x(0)=x_0 ,
\end{equation}
with the multivalued mapping $\mathcal{K}f:I\times\mathbb{R}^n\rightarrow \mathcal{P}(\mathbb{R}^n)$ defined as
\begin{equation}\label{eq_Kf}
\mathcal{K}f(t,x)=\bigcap_{\varepsilon>0}\overline{\rm co}f\left(t,\overline{B}_{\varepsilon}(x)\right) \quad \text{ for every } (t,x)\in I\times\mathbb{R}^{n},
\end{equation}
where $\overline{\rm co}$ means closed convex hull, and $\overline{B}_{\varepsilon}(x)$ is the closed ball centered at $x$ and radius $\varepsilon>0$.

Absolutely continuous solutions of \eqref{ivpI} are usually called \textit{Krasovskij solutions} of \eqref{ivp}. We need the following known result concerning the existence of Krasovskij solutions, see \cite[Proposition 2.1]{CidPo}. 

\begin{proposition}\label{prop_Kras}
	Let $f:I\times\mathbb{R}^n\rightarrow\mathbb{R}^n$ be such that
	\begin{enumerate}
		\item[$(1)$] There exists $M\in L^1(I)$ such that for a.a. $t\in I$ and all $x\in \mathbb{R}^n$, we have
		\[\left\|f(t,x)\right\|\leq M(t)(1+\left\|x\right\|); \]
		\item[$(2)$] For all $x\in \mathbb{R}^n$, $f(\cdot,x)$ is measurable.
	\end{enumerate}
	Then problem \eqref{ivp} has at least one Krasovskij solution.
\end{proposition}

\section{Existence result}\label{sec_exis}

Our aim is to prove the existence of solutions to \eqref{ivp} and to localize them in a viable region $R$. In order to do that, we introduce the following notion which links the differential problem \eqref{ivp} and the viable region.

\begin{definition}\label{def_solR}
	A set $R\subset I\times\R^n$ is called a \textit{solution-region} of \eqref{ivp} if it is a viable region with the associated viable pair $(h,p)$ satisfying the following conditions:
	\begin{enumerate}
		\item[$(i)$] for a.a. $t$ and all $x$ with $(t,x)\notin R$, one has
		\begin{equation}\label{eq_solR}
		\langle \nabla h(t,x), (1,z)\rangle\leq 0, \quad \text{for all } z\in\mathcal{K}f(p(t,x));
		\end{equation} 
		\item[$(ii)$] $h(0,x_0)\leq 0$.
	\end{enumerate}
\end{definition}

\begin{remark}
In case that $f(t,\cdot)$ be continuous at $x$, then $\mathcal{K}f(t,x)=\{f(t,x) \}$, see \cite{aubin}. Thus, if $f$ is continuous with respect to the second variable, inequality \eqref{eq_solR} is equivalent to
\begin{equation}\label{eq_solR_cont}
\langle \nabla h(t,x), (1,f(p(t,x)))\rangle\leq 0,
\end{equation}	
as required in \cite{f,ftt,t} for the definition of a solution-region of \eqref{ivp}.

Condition \eqref{eq_solR_cont} resembles the one in the notion of a \textit{bound set} relative to \eqref{ivp} due to Gaines and Mawhin \cite{gm}. Note that, unlike \cite{gm}, here $h$ needs not be negative in the interior of the set $K$. In fact, $K$ may have empty interior. 
\end{remark}	

The following technical result concerning measure theory is crucial. It is a straightforward consequence of \cite[Lemma 5.8.13]{bo}.

\begin{lemma}
	\label{lebo1}
	Let $a, \, b \in \mathbb R$, $a<b$. 
	
	If $\varphi:[a,b] \longrightarrow \mathbb R$ is almost everywhere differentiable on $[a,b]$, then for any null measure set $A \subset \mathbb R$ there exists a null measure set $B \subset \varphi^{-1}(A)$ such that
	$$\varphi'(t)=0 \quad \mbox{for all $t \in \varphi^{-1}(A) \setminus B$.}$$
\end{lemma}

Now we are in a position to establish our main existence and localization result for \eqref{ivp}.

\begin{theorem}\label{th}
	Assume that there exists $R$ a solution-region of \eqref{ivp}.
	
	Let $f:I\times\R^n\rightarrow\R^n$ be such that
	\begin{enumerate}
		\item[$(H_1)$] for each $\mu\geq 0$, there exists $M_{\mu}\in L^1(I)$ such that for a.a. $t\in I$ and all $x\in\R^n$ with $\left\|x\right\|\leq\mu$, one has
		\[\left\|f(t,x)\right\|\leq M_{\mu}(t); \] 
		\item[$(H_2)$] for all $x\in\R^n$, $f(p(\cdot,x))$ is measurable;
		\item[$(H_3)$] there exists a countable number of continuously differentiable functions $\tau_n:I\times\R^n\rightarrow\R$ and a countable number of null-measure sets $A_n\subset\R$ $(n\in\N)$ such that for a.a. $t\in I$,
		\[f(t,x) \text{ is continuous in } x \text{ for each } (t,x)\in R\setminus\bigcup_{n=1}^{\infty}\tau_n^{-1}(A_n) \]
		and for each $n\in\N$ and for each $(t,x)\in\tau_n^{-1}(A_n)\cap R$, we have
		\begin{equation}\label{eq_transv}
		\langle \nabla\tau_n(t,x),(1,z)\rangle\neq 0 \quad \text{for all } z\in\mathcal{K}f(t,x).
		\end{equation}
	\end{enumerate}
	Then problem \eqref{ivp} has at least one absolutely continuous solution $x$ such that $(t,x(t))\in R$ for every $t\in I$.
\end{theorem}

\noindent
{\bf Proof.} Consider the auxiliary problem
\begin{equation}
\label{ivp_modI}
x'(t)\in\mathcal{K}\tilde{f}(t,x(t)) \quad \text{for a.a. } t \in I, \ \ x(0)=x_0 ,
\end{equation}
where $\tilde{f}(t,x)=f(p(t,x))$ for every $(t,x)\in I\times\R^n$.

Since $p$ is bounded and $f$ satisfies condition $(H_1)$, there exists $\tilde{M}\in L^1(I)$ such that for a.a. $t\in I$ and all $x\in \R^n$, one has
\[\left\|\tilde{f}(t,x)\right\|\leq \tilde{M}(t). \]
In addition, assumption $(H_2)$ implies that for all $x\in\R^n$, $\tilde{f}(\cdot,x)$ is measurable. Hence, Proposition \ref{prop_Kras} ensures that the modified problem \eqref{ivp_modI} has at least one absolutely continuous solution.

Now, let us divide the proof into two steps.

\medskip

\textbf{Step 1.} If $x$ is a solution of problem \eqref{ivp_modI}, then $(t,x(t))\in R$ for every $t\in I$.

Assume the contrary. Then, since by Definition \ref{def_solR} $(ii)$, $h(0,x(0))=h(0,x_0)\leq 0$, there exist $t_1,t_2\in[0,T]$, $t_1<t_2$, such that $h(t_1,x(t_1))= 0$ and $h(t,x(t))>0$ for all $t\in(t_1,t_2)$. By the chain rule, it follows that almost everywhere on $(t_1,t_2)$, one has
\begin{equation*}
\frac{dh}{dt}(t,x(t))= \langle\nabla h(t,x(t)), (1, x'(t))\rangle =\langle\nabla h(t,x(t)), (1, z)\rangle,
\end{equation*} 
where $z\in\mathcal{K}\tilde{f}(t,x(t))$. Since $\mathcal{K}\tilde{f}(t,x(t))\subset \mathcal{K}f(p(t,x(t)))$, it follows from condition $(i)$ in Definition \ref{def_solR} that almost everywhere on $(t_1,t_2)$,
\[\frac{dh}{dt}(t,x(t))\leq 0. \]
Therefore, $h(t,x(t))\leq 0$ for all $t\in [t_1,t_2]$, a contradiction. Then $(t,x(t))\in R$ for all $t\in I$. 

\medskip

\textbf{Step 2.} All the solutions of problem \eqref{ivp_modI} are solutions of the initial value problem \eqref{ivp}.

By Step 1, solutions of \eqref{ivp_modI} solve
\[x'(t)\in\mathcal{K}f(t,x(t)) \quad \text{for a.a. } t \in I, \ \ x(0)=x_0 , \]
since they are located in $R$, $p(t,x)=(t,x)$ for every $(t,x)\in R$ and $\mathcal{K}\tilde{f}(t,x)\subset\mathcal{K}f(p(t,x))$ for every $(t,x)\in I\times\R^n$.

Let $x$ be a solution of problem \eqref{ivp_modI}. Let us prove that for each $n\in\N$,
\[m\left(\{t\in I:\tau_n(t,x(t))\in A_n \} \right)=0,\]
(where $m$ stands for the Lebesgue measure). Denote $J_n=\{t\in I:\tau_n(t,x(t))\in A_n \}$ and $\varphi_n(t)=\tau_n(t,x(t))$ for all $t\in I$. Since $m(A_n)=0$, by Lemma \ref{lebo1}, there exists $B_n\subset\varphi_n^{-1}(A_n)$, with $m(B_n)=0$, such that for every $t\in\varphi^{-1}(A_n)\setminus B_n$, one has
\[\frac{d\tau_n(t,x(t))}{dt}=0. \]
Thus, 
\[\langle\nabla \tau_n(t,x(t)) ,(1,x'(t))\rangle=0 \quad \text{ for all } t\in \varphi_i^{-1}(A_n) \setminus B_n. \]
Now, condition \eqref{eq_transv} implies that $\varphi_i^{-1}(A_n)$ is a null-measure set.

Therefore, if $x$ is a solution of \eqref{ivp_modI}, then for a.a. $t\in I$ the mapping $f(t,\cdot)$ is continuous at $x(t)$, which implies that $\mathcal{K}f(t,x(t))=\{f(t,x(t)) \}$ for a.a. $t\in I$. Then, $x'(t)\in \mathcal{K}f(t,x(t))=\{f(t,x(t)) \}$ for a.a. $t\in I$, so $x$ is an absolutely continuous solution of \eqref{ivp}.
\qed

\begin{remark}\label{rmk_def_solR}
	Observe that the concept of solution-region of problem \eqref{ivp} is just needed in the Step 1 of the proof of Theorem \ref{th}. In view of this part of the proof, it is clear that condition $(i)$ in Definition \ref{def_solR} can be replaced by the following one:
	 \begin{enumerate}
	 	\item[$(i^*)$] for a.a. $t$ and all $x$ with $(t,x)\notin R$, one has
	 	\begin{equation*}
	 	\langle\nabla h(t,x), (1,z)\rangle\leq 0, \quad \text{for all } z\in\mathcal{K}\tilde{f}(t,x).
	 	\end{equation*} 
	 \end{enumerate}
\end{remark}

\begin{remark}
	Notice that condition $(H_3)$ in Theorem \ref{th} is satisfied in particular if the function $f$ can be expressed in the form
	\begin{equation}
	\label{efe}
	f(t,x)=F(t,g_1(\tau_1(t,x),x),g_2(\tau_2(t,x),x),\dots,g_N(\tau_N(t,x),x))  \quad \mbox{for some $N \in \mathbb N$,}
	\end{equation}
	where:
	\begin{itemize}
			\item each function $g_i:\mathbb R \times \mathbb R^n \longrightarrow \mathbb R$ is continuous in $(\mathbb{R}\setminus A_i)\times\mathbb{R}^n$;
			\item for a.a. $t\in I$, $F(t,\cdot)$ is continuous in $\mathbb R^N$.
	\end{itemize}
	The particular form of $f$ given by \eqref{efe} was considered in \cite{pmr,pr} and, previously, in \cite{BreSh}. Moreover, the transversality condition \eqref{eq_transv} is a localized version of that due to Bressan and Shen \cite{BreSh}. Similar conditions can also be found in \cite{AGG,AJ} in case of discontinuous differential inclusions and periodic problems, respectively.
\end{remark}

	Both conditions given by inequalities \eqref{eq_solR} and \eqref{eq_transv} can be thought as viability conditions for the solutions of \eqref{ivp}. Indeed, \eqref{eq_solR} is the key ingredient in order to ensure that the graphs of the solutions of the modified problem \eqref{ivp_modI} are located in $R$ (see Step 1 of the proof of Theorem \ref{th}), whereas \eqref{eq_transv} is crucial to guarantee that the solutions of \eqref{ivp_modI} avoid the hyper-surfaces of type $\tau_n(t,x)=c$ with $c\in A_n$, where the function $f$ may be discontinuous (see Step 2 of the proof of Theorem \ref{th}). 

\section{Examples}\label{sec_examples}

\subsection{Localizing the solutions in a ball}

Let $R\subset I\times\R^n$ be such that $R=I\times \overline{B}_r$, where $\overline{B}_r$ is the closed ball centered at the origin with radius $r>0$. Observe that $R$ is clearly a viable region with the viable pair $(h,p)$ given by
\begin{itemize}
	\item $p(t,x)=(t,p_r(x))$, where $p_r:\R^n\rightarrow \overline{B}_r$ is the orthogonal projection onto the ball $\overline{B}_r$, i.e.,
	\[p_r(x)=\left\{\begin{array}{ll} x, & \text{ if }x\in \overline{B}_r, \\[0.2cm] r\displaystyle\frac{x}{\left\|x\right\|}, & \text{ if } x\notin \overline{B}_r. \end{array} \right. \]
	\item $h(t,x)=h(x)=\displaystyle\frac{1}{2}\left\|x-p_r(x)\right\|^2$.
\end{itemize}
Note that $h$ is a continuously differentiable function and $\nabla h(x)=x-p_r(x)$ for all $x\in\R^n$, see \cite{t}.

Let $x_0\in \overline{B}_r$. The set $R$ will be a solution region of \eqref{ivp} with the pair $(h,p)$ provided that for a.a. $t\in I$ and all $x\in\R^n\setminus \overline{B}_r$, one has
\begin{equation*}
\langle x-p_r(x), z\rangle\leq 0, \quad \text{for all } z\in\mathcal{K}f(t,p_r(x)),
\end{equation*} 
or equivalently,
\begin{equation*}
\langle r\frac{x}{\left\|x\right\|},z\rangle\leq 0, \quad \text{for all } z\in\mathcal{K}f\left(t,r\frac{x}{\left\|x\right\|}\right).
\end{equation*} 
In other words, $R$ will be a solution region of \eqref{ivp} with the pair $(h,p)$ provided that for a.a. $t\in I$ and all $x\in\R^n$ with $\left\|x\right\|=r$, one has
\begin{equation}\label{eq_h_ball}
\langle x, z\rangle\leq 0, \quad \text{for all } z\in\mathcal{K}f\left(t,x\right).
\end{equation}
Therefore, inequality \eqref{eq_h_ball} gives us a simple sufficient condition in order to ensure that the set $I\times \overline{B}_r$ is a solution region of the initial value problem \eqref{ivp}.

Now we illustrate the existence result with the following example where the solution is localized in the unit ball.

\begin{example}
Consider the initial value problem 
\begin{equation}\label{eq_exBall}
\left\{\begin{array}{lll} x'=x^3+y-3\,x+\phi(x^2+y^2+\alpha\, t), & \quad t\in[0,1], & \quad x(0)=0, \\ y'=y^3-x-3\,y\,e^{|x|}, & \quad t\in[0,1], & \quad y(0)=0, \end{array} \right.
\end{equation}	
where $\alpha\in\R$ and the function $\phi:\R\rightarrow\R$ is continuous at irrational numbers and discontinuous at every rational number, and such that $0<\phi(s)<1$ for all $s\in\R$.

We shall prove that \eqref{eq_exBall} has at least one absolutely continuous solution located in the ball $\overline{B}_1$ provided that $|\alpha|$ is sufficiently large.

Denote $g(x,y)=x^3+y-3\,x$ and $f=(f_1,f_2)$, with $f_1(t,x,y)=g(x,y)+\phi(x^2+y^2+\alpha\, t)$ and $f_2(t,x,y)=y^3-x-3\,y\,e^{|x|}$. 

First, let us show that the set $[0,1]\times\overline{B}_1$ is a solution region of the initial value problem \eqref{eq_exBall} with the viable pair $(h,p)$ given above. 
Let $(x,y)\in\R^2$ with $\left\|(x,y)\right\|=1$. Since $g$ and $f_2$ are continuous functions, there exists $\varepsilon>0$ such that
\[x\,g_1(\bar{x},\bar{y})+y\,f_2(\bar{x},\bar{y})\leq x\,g_1(x,y)+y\,f_2(x,y)+1 \quad \text{for every } (\bar{x},\bar{y})\in B_{\varepsilon}((x,y)), \]
and so
\[x\,g_1(\bar{x},\bar{y})+y\,f_2(\bar{x},\bar{y})\leq x^4-3\,x^2+y^4-3\,y^2\,e^{|x|}+1\leq (x^2+y^2)^2-3(x^2+y^2)+1=-1, \]
for every $(\bar{x},\bar{y})\in B_{\varepsilon}((x,y))$.
Now, since $0<\phi(s)<1$ for all $s\in\R$, 
\[x\,f_1(t,\bar{x},\bar{y})+y\,f_2(\bar{x},\bar{y})\leq-1+x\,\phi(\bar{x}^2+\bar{y}^2+\alpha\, t)\leq 0 \quad \text{for every } (\bar{x},\bar{y})\in B_{\varepsilon}((x,y)). \]
Hence, one has that 
\begin{equation}\label{eq_ex_h_ball}
\langle (x,y), (z_1,z_2)\rangle\leq 0 \quad \text{for all } z=(z_1,z_2)\in f(t,\overline{B}_{\varepsilon}((x,y))).
\end{equation}
Observe that \eqref{eq_ex_h_ball} implies that 
\[\langle (x,y), z\rangle\leq 0 \quad \text{for all } z\in \mathcal{K}f\left(t,(x,y)\right)\subset \overline{\rm co}f(t,\overline{B}_{\varepsilon}((x,y))).  \]
Indeed, if \eqref{eq_ex_h_ball} holds, then for $z,\bar{z}\in f(t,\overline{B}_{\varepsilon}((x,y)))$ and $\lambda\in[0,1]$, one has
\[\langle (x,y), \lambda\, z+(1-\lambda)\bar{z}\rangle=\lambda\langle(x,y), z\rangle+(1-\lambda)\langle(x,y),\bar{z}\rangle\leq 0. \]
Then the set $R=[0,1]\times\overline{B}_1$ is a solution region of the initial value problem \eqref{eq_exBall}.

Finally, the existence of an absolutely continuous solution for problem \eqref{eq_exBall} with a large $|\alpha|$ is obtained as a consequence of Theorem \ref{th}. Clearly, $f$ satisfies assumptions $(H_1)$ and $(H_2)$. For condition $(H_3)$, choose the function $\tau(t,x,y)=x^2+y^2+\alpha\, t$ and as $A$ the set of all rational numbers. By the definition of $\phi$, $f$ is continuous in $R\setminus\tau^{-1}(A)$. Since $f$ is bounded in $R$, so is $\mathcal{K}f$ and thus there exists $M>0$ such that 
\[|x\,z_1+y\,z_2|\leq M \quad \text{for all } (z_1,z_2)\in \mathcal{K}f\left(t,(x,y)\right) \text{ and every } (t,x,y)\in R.\]
Therefore, if $|\alpha|>2\,M$, then for every $(t,x,y)\in R$, we have
\[\langle\nabla\tau(t,x,y),(1,z_1,z_2)\rangle=\alpha +2(x\,z_1+y\,z_2)\neq 0 \quad \text{for all } (z_1,z_2)\in \mathcal{K}f\left(t,(x,y)\right). \]
In conclusion, Theorem \ref{th} ensures that problem \eqref{eq_exBall} has at least one absolutely continuous solution $(x,y)$ such that $(x(t),y(t))\in \overline{B}_1$ for all $t\in[0,1]$ provided that $|\alpha|>2\,M$.
\end{example}

Obviously, the choice of the function $h$ when it comes to defining the viable pair $(h,p)$ associated to $R$ is not unique. Distinct choices of $h$ derive into different conditions on $f$ in order to check that $R$ is a solution region of \eqref{ivp}. This idea allows us to state the following existence result.

\begin{proposition}
	Assume that there exist $r>0$ and a continuously differentiable function $h:\R^n\rightarrow\R$ such that $\overline{B}_r=h^{-1}((-\infty,0])$ and for a.a. $t\in I$ and all $x\in\R^n$ with $\left\|x\right\|=r$, one has
	\begin{equation*}
	\langle \nabla h(\lambda\,x), z\rangle\leq 0, \quad \text{for all } z\in\mathcal{K}f\left(t,x\right) \text{ and all } \lambda\geq 1.
	\end{equation*}
	Moreover, assume that $f:I\times\R^n\rightarrow\R^n$ satisfies conditions $(H_1)$, $(H_2)$ and $(H_3)$. 
	
	Then the initial value problem \eqref{ivp} has at least one absolutely continuous solution $x$ such that $\left\|x(t)\right\|\leq r$ for all $t\in I$ provided that $x_0\in\overline{B}_r$.
\end{proposition}

To the best of our knowledge, the previous result is even new in the classical case of Carath\'eodory nonlinearities.

\begin{corollary}
	Let $f:I\times\R^n\rightarrow\R^n$ be a Carath\'eodory function. 
	
	Assume that there exist $r>0$ and a continuously differentiable function $h:\R^n\rightarrow\R$ such that $\overline{B}_r=h^{-1}((-\infty,0])$ and for a.a. $t\in I$ and all $x\in\R^n$ with $\left\|x\right\|=r$, one has
	\begin{equation*}
	\langle \nabla h(\lambda\,x), f(t,x)\rangle\leq 0, \quad \text{for all } \lambda\geq 1.
	\end{equation*}	
	
	Then the initial value problem \eqref{ivp} has at least one absolutely continuous solution $x$ such that $\left\|x(t)\right\|\leq r$ for all $t\in I$ provided that $x_0\in\overline{B}_r$.
\end{corollary}

\subsection{Discontinuous lower and upper solutions}

First, let us introduce the notion of absolutely continuous lower and upper solutions for problem \eqref{ivp} in the scalar case (i.e., $n=1$).

\begin{definition}\label{def_lu}
	A function $\alpha\in W^{1,1}(I,\R)$ is called a lower solution of \eqref{ivp} if $\alpha'(t)\leq f(t,\alpha(t))$ for a.e. $t\in I$ and $\alpha(0)\leq x_0$.
	An upper solution of \eqref{ivp} is defined analogously reversing the previous inequalities.
\end{definition}

Assume that $\alpha,\beta\in W^{1,1}(I,\R)$ are, respectively, a lower and an upper solution of \eqref{ivp} with $\alpha\leq \beta$. Then the region
\[R=\{(t,x)\in I\times\R:\alpha(t)\leq x\leq \beta(t) \} \]
is a solution region of \eqref{ivp}. Indeed, the pair $(h,p)$ given by
\[h(t,x)=\max\{x-\beta(t),\alpha(t)-x,0 \}, \qquad p(t,x)=(t,\max\{\min\{x,\beta(t) \},\alpha(t) \}), \]
is a viable pair satisfying that $h(0,x_0)\leq 0$ and condition $(i^*)$ in Remark \ref{rmk_def_solR}. Hence, Theorem 2.5 in \cite{pmr} can be seen as a consequence of Theorem \ref{th}. In the case of a Carath\'eodory function $f$, the details can be found in \cite{ftt}.

The problem of finding a lower and an upper solution for \eqref{ivp} is not easy in general and so in most of cases limit us to consider constant or linear functions. In this way, it could be convenient to relax the regularity assumptions in Definition \ref{def_lu} in order to consider piecewise constant or linear functions. Piecewise continuous lower and upper solutions can be found in a few papers in the literature, see \cite{bp,lp,top}.

The next example shows that in some cases the region between well--ordered discontinuous lower and upper solutions also becomes a solution region, which provides existence directly from Theorem \ref{th}. It was inspired by \cite[Example 4.9]{top}.

\begin{example}\label{ex_LU}
	Consider the initial value problem 
	\begin{equation}\label{eq_exLU}
	x'=-x^2-x+2\,t+1, \quad t\in[0,1], \quad x(0)=0. 
	\end{equation}
	
	The set $R=\{(t,x)\in [0,1]\times\R:\alpha(t)\leq x\leq \beta(t) \}$, with $\alpha,\beta:[0,1]\rightarrow\R$ given by
	\[\alpha(t)=t, \qquad \beta(t)=\left\{\begin{array}{ll} 1, & \text{ if } t\in[0,1/2), \\ 2, & \text{ if } t\in[1/2,1], \end{array} \right. \]
	is a solution region of \eqref{eq_exLU}. Indeed, let us define the maps $h:[0,1]\times\R\rightarrow\R$ and $p:[0,1]\times\R\rightarrow [0,1]\times\R$ by
	\[h(t,x)=\left\{\begin{array}{ll} 0, & \text{ if } (t,x)\in R, \\ t-x, & \text{ if } x<t, \\ (x-1)(1/2-t), & \text{ if } 0\leq t<1/2, \ 1<x\leq 2\,t+1, \\ (x-1)(1/2-t)+(x-2\,t-1)^2, & \text{ if } 0\leq t<1/2, \ 2\,t+1<x, \\ (x-2)^2, & \text{ if } 1/2\leq t\leq 1, \ 2<x, \end{array} \right. \]
	and
	\[p(t,x)=\left\{\begin{array}{ll} (t,t), & \text{ if } x<t, \\ (t,x), & \text{ if } t\leq x\leq \min\{2\,t+1,2 \}, \\ (t,\min\{2\,t+1,2 \}), & \text{ if } x>\min\{2\,t+1,2 \}. \end{array} \right. \]
	The map $p$ is continuous and such that $p(t,x)=(t,x)$ for every $(t,x)\in R$, and $h$ is a nonnegative continuous function with $h^{-1}(0)=R$. Also,  $\frac{\partial h}{\partial t}$ and $\frac{\partial h}{\partial x}$ are Carath\'eodory maps on $R^c$. Hence, $(h,p)$ is a viable pair associated to $R$. Moreover, one may easily check that for a.a. $t$ and all $x$ with $(t,x)\in R^c$,
	\[\frac{\partial h}{\partial t}(t,x)+\frac{\partial h}{\partial x}(t,x)f(p(t,x))\leq 0, \]
	where $f(t,x)=-x^2-x+2\,t+1$. Therefore, $R$ is a solution region of \eqref{eq_exLU} and, since $f$ is continuous, Theorem \ref{th} immediately ensures the existence of a solution for \eqref{eq_exLU} between $\alpha$ and $\beta$. 
	
	Observe that $\beta$ cannot be an upper solution in the sense of Definition \ref{def_lu} since it is not continuous. Note that the constant function $\gamma(t)=1$, $t\in [0,1]$, is not an upper solution of \eqref{eq_exLU}. Even if the absolutely continuous function defined as $\tilde{\beta}(t)=\min\{2\,t+1,2 \}$ is clearly an upper solution of \eqref{eq_exLU}, $\beta$ provides a sharper localization of the solution than $\tilde{\beta}$. 
	
	On the other hand, $(h,p)$ is not an admissible pair associated to $R$ since assumption $(iii)$ in Definition \ref{def_adm} does not hold. Indeed, $p_2(t,x)=x$ on the set $\{(t,x)\in [0,1/2)\times\R: 1<x\leq 2\,t+1 \}\subset R^c$, which implies
	\[\langle \nabla_x h(t,x),p_2(t,x)-x \rangle=\frac{\partial h}{\partial x}(t,x)\left(p_2(t,x)-x \right)=0. \]
	Therefore, the theory developed in \cite{f} cannot be applied here with this pair $(h,p)$.
\end{example}

Finally, note that the existence of well--ordered piecewise continuous lower and upper solutions, with downwards (upwards) jumps in case of the lower (upper) solutions, implies the existence of a better pair of well--ordered absolutely continuous lower and upper solutions for \eqref{ivp}, as proven in \cite[Proposition 4.2]{bp}. Therefore, the solution-region can be defined by means of the absolutely continuous lower and upper solutions. In this way, Theorem 2.1 in \cite{lp} can also be seen as a corollary of the main result in this paper.

\section*{Acknowledgements}

The author is grateful to Prof. Rodrigo L\'opez Pouso for his steady encouragement and helpful suggestions.

He was partially supported by Xunta de Galicia (Spain), project ED431C 2019/02 and AIE, Spain and FEDER, grant PID2020-113275GB-I00.

\end{document}